\documentclass[12pt,leqno,a4paper]{article}
\usepackage{amsmath,amsthm}
\usepackage[latin1]{inputenc}
\usepackage[T1]{fontenc}
\usepackage[english]{babel}

\newtheorem{theorem}{Theorem}[section]
\newtheorem{proposition}[theorem]{Proposition}
\newtheorem{lemma}[theorem]{Lemma}
\newtheorem{corollary}[theorem]{Corollary}



\newcommand*\proofnamestyle{\itshape}

\newcounter{std}

\begin{document}

    \title{Differential analysis of matrix convex functions}
      \author{Frank Hansen and Jun Tomiyama}
      \date{January 12, 2006\\
            {\tiny Revised June 12, 2006}}
      \maketitle

      \begin{abstract}
      We analyze matrix convex functions of a fixed order defined in a real interval
    by differential methods as opposed to the characterization in terms
    of divided differences given by Kraus~\cite{kn:kraus:1936}. We obtain for each order conditions for matrix
    convexity which are necessary and locally sufficient, and they
    allow us to prove the existence of gaps between classes of matrix convex functions of
    successive orders, and to give explicit examples of the type of functions contained in each of
    these gaps. The given conditions are shown to be also globally sufficient for matrix convexity of order two.
    We finally introduce a fractional transformation which connects the set of matrix monotone functions
    of each order $ n $ with the set of matrix convex functions of the following order $ n+1. $
    \end{abstract}

    \section{Introduction}

    Let $ f $ be a real function defined in an interval $ I. $ It is said to be
    $ n $-convex if
    \[
    f(\lambda A+(1-\lambda)B)\le \lambda f(A)+(1-\lambda)f(B)\qquad\lambda\in[0,1]
    \]
    for arbitrary Hermitian $ n\times n $ matrices $ A $ and $ B $ with spectra in $ I. $ It is said to be
    $ n $-concave if $ -f $ is $ n $-convex, and it is said to be $ n $-monotone if
    \[
    A\le B\quad\Rightarrow\quad f(A)\le f(B)
    \]
    for arbitrary Hermitian $ n\times n $ matrices $ A $ and $ B $ with spectra in $ I. $
    We consider the interval of definition to be part of the specification of a function.
    The notions of $ n $-convexity and $ n $-monotonicity are therefore not associated solely with an assignment rule,
    but may depend on the interval in which the rule is applied. We denote by $ P_n(I) $ the set
    of $ n $-monotone functions defined in an interval $ I, $ and by $ K_n(I) $ the set of $ n $-convex functions
    defined in $ I. $

    We shall sometimes use the standard regularization procedure, cf. for example
    Donoghue \cite[Page 11]{kn:donoghue:1974}.
    Let $ \varphi $ be a positive and even $ C^{\infty} $-function defined on the real axis, vanishing
    outside the closed interval $ [-1,1] $ and normalized such that
    \[
    \int_{-1}^1 \varphi(x)\,dx=1.
    \]
    For any locally integrable function $ f $ defined in an open interval $ (a,b) $ we form its regularization
    \[
    f_\epsilon(t)=\frac{1}{\epsilon}\int_a^b \varphi\left(\frac{t-s}{\epsilon}\right)f(s)\,ds\qquad t\in{\mathbf R}
    \]
    for small $ \epsilon>0, $ and realize that it is infinitely many times differentiable.
    For $ t\in(a+\epsilon,b-\epsilon) $ we may also write
    \[
    f_\epsilon(t)=\int_{-1}^1\varphi(s)f(t-\epsilon s)\,ds.
    \]
    If $ f $ is continuous, then $ f_\epsilon $ converges uniformly to $ f $ on any compact subinterval
    of $ (a,b). $ If in addition $ f $ is $ n $-convex (or $ n $-monotone) in $ (a,b), $
    then $ f_\epsilon $ is $ n $-convex (or $ n $-monotone) in the slightly smaller interval
    $ (a+\epsilon,b-\epsilon). $ Since the pointwise
    limit of a sequence of $ n $-convex (or $ n $-monotone) functions is again $ n $-convex (or $ n $-monotone),
    we may therefore in many applications assume that an $ n $-convex or $ n $-monotone function is
    sufficiently many times differentiable.

    We remind the reader \cite[Chapter VII Theorem VI]{kn:donoghue:1974} that a $ 2 $-monotone function
    defined in an open interval automatically is continuously differentiable. The first statement of the
    following lemma can also be found in the same reference.

    \begin{lemma}\label{lemma: 2-monotone and 2-convex functions}
    Let $ I $ be an open interval and consider a real function $ f $ defined in $ I. $
    \begin{list}{(\arabic{std})}{\usecounter{std}}

    \item If $ f $ is $ 2 $-monotone and $ f'(t_0)=0 $ for a single $ t_0\in I, $ then $ f $ is a constant
    and therefore $ f'(t)=0 $ for all $ t\in I. $

    \item If $ f $ is twice continuously differentiable, $ 2 $-convex and $ f''(t_0)=0 $ for a single
    $ t_0\in I, $ then $ f $ is affine and therefore $ f''(t)=0 $ for all $ t\in I. $

    \end{list}
    \end{lemma}

    \begin{proof}
    We use the characterizations of $ 2 $-monotonicity by Löwner \cite{kn:loewner:1934} and
    $ 2 $-convexity by Kraus \cite{kn:kraus:1936}. If $ f $ is $ 2 $-monotone then
    \[
    \det\begin{pmatrix}
        [t_0,t_0]_f & [t_0,t]_f\\
        [t,t_0]_f   & [t,t]_f
        \end{pmatrix}=f'(t_0)f'(t)-[t_0,t]_f^2\ge 0
    \]
    from which the first statement follows. Similarly, if $ f $ is $ 2 $-convex then
    \[
    \det\begin{pmatrix}
        [t_0,t_0,t_0]_f & [t_0,t_0,t]_f\\
        [t_0,t,t_0]_f   & [t_0,t,t]_f
        \end{pmatrix}=\frac{1}{2}f''(t_0)[t_0,t,t]_f-[t_0,t_0,t]_f^2\ge 0
    \]
    from which the second statement follows.

    \end{proof}

    \subsection{Differential conditions}

    If $ I $ is open and $ f $ is twice differentiable
    we introduce for real numbers $ s,t_1,\dots,t_n $ in $ I $ the leading determinants
    \begin{gather}\label{leading determinants}
    D_r(s)=\det
    \left(\begin{array}{cccc}
        [t_1, s, t_1]_f & [t_1, s, t_2]_f & \cdots & [t_1, s, t_r]_f\\[1ex]
        [t_2, s, t_1]_f & [t_2, s, t_2]_f & \cdots & [t_2, s, t_r]_f\\[1ex]
        \vdots      & \vdots        &     & \vdots\\[1ex]
        [t_r, s, t_1]_f & [t_r, s, t_2]_f & \cdots & [t_r, s, t_r]_f
        \end{array}\right)
    \end{gather}
    for $ r=1,\dots,n. $ Kraus proved \cite{kn:kraus:1936} that $ f $
    is $ n $-convex if and only if the matrix
    \begin{equation}\label{Kraus matrix}
    H(s)=\Big([t_i, s, t_j]_f\Big)_{i,j=1}^n
    \end{equation}
    is positive semi-definite for every sequence $ t_1,\dots,t_n\in I $ and
    $ s=t_1,\dots,t_n. $

    \begin{theorem}\label{determinant condition}
    Let $ f $ be a real $ 2n $ times continuously differentiable function defined in an open interval
    $ I. $ The determinants defined in (\ref{leading determinants}) can be written
    \[\label{leading determinants reformulated}
    D_r(s)=\det M_r\cdot\prod_{k=1}^{r-1}\,\,\prod_{l=1}^{r-k}(t_{k+l}-t_l)^2\qquad r=1,\dots,n
    \]
    where $ M_r=\big(m_{ij}\bigr)_{i,j=1}^r $ and $  m_{ij}=[t_1,\dots, t_i, s, t_1,\dots, t_j]_f. $
    If in addition $ f $ is $ n $-convex, then the matrix
    \[
    K_n(f;t)=\begin{pmatrix}
    \displaystyle\frac{f^{i+j}(t)}{(i+j)!}
    \end{pmatrix}_{i,j=1}^n
    \]
    is positive semi-definite for each $ t\in I. $ On the other hand, if
    $ K_n(f;t_0) $ is positive definite for some $ t_0\in I, $ then $ f $ is $ n $-convex
    in some open interval $ J $ with $ t_0\in J\subseteq I. $
    \end{theorem}

    \begin{proof}
    Let $ t_1,\dots,t_n $ be $ n $ distinct points in the interval $ I. $
    Since there is no possibility of confusion, we shall omit the reference to the function $ f $
    in the divided differences. For each $ r=2,\dots,n $ we intend to prove
    \begin{gather}\label{intermediate determinant identity}
    D_r(s)=\det M_r(p)\cdot\prod_{k=1}^p\,\,\prod_{l=1}^{r-k}(t_{k+l}-t_l)^2
    \end{gather}
    by induction for $ p=1,\dots, r-1 $ where
    \[
    M_r(p)=\left(\begin{array}{cc}
           M_p & \bigl([t_1,\dots,t_i, s, t_{j-p},\dots,t_j]\bigr)_{i,j}\\[1ex]
           \bigl([t_{i-p},\dots,t_i, s, t_1,\dots,t_j]\bigr)_{i,j}
           & \bigl([t_{i-p},\dots,t_i, s, t_{j-p},\dots,t_j]\bigr)_{i,j}
           \end{array}\right).
    \]
    Note that $ M_r(p) $ is an $ r\times r $ matrix written as a block matrix with the
    $ p\times p $ matrix $ M_p $ as the $ (1,1) $ block entry. The indices $ i $ and $ j $ refer to
    the absolute row and column numbers in $ M_r(p). $ The row index $ i $ in block entry $ (2,1) $ hence runs from
    $ p+1 $ to $ r, $ and the column index $ j $ runs from 1 to $ p. $ When
    (\ref{intermediate determinant identity}) is proved, the
    first part of the theorem follows by setting $ p=r-1 $ and noting that $ M_r(r-1)=M_r. $

    In the determinant expression (\ref{leading determinants}) we subtract the first row from the second row,
    the second row from the third and so forth until the $ (r-1) $th row is subtracted from the $ r $th row.
    We thus obtain
    \[
    \begin{array}{l}
    D_r(s)=\\[2ex]
    \left|\begin{array}{cccc}
[t_1, s, t_1]                   & [t_1, s, t_2]                   & \cdots & [t_1, s, t_r]\\[1ex]
[t_2, s, t_1]-[t_1, s, t_1]     & [t_2, s, t_2]-[t_1, s, t_2]     & \cdots & [t_2, s, t_r]-[t_1, s, t_r]\\[1ex]
        \vdots      & \vdots                      &        & \vdots\\[1ex]
[t_r, s, t_1]-[t_{r-1}, s, t_1 ]& [t_r, s, t_2]-[t_{r-1}, s, t_2] & \cdots & [t_r, s, t_r]-[t_{r-1}, s, t_r]
    \end{array}\right|
    \end{array}
    \]
    and since for $ i=2,\dots,r $ and $ j=1,\dots,r $ the difference
    \[
    \begin{array}{rl}
    [t_i, s, t_j]-[t_{i-1}, s, t_j]&=[t_i, s, t_j]-[s, t_j, t_{i-1}]
    =(t_i-t_{i-1})[t_i, s, t_j, t_{i-1}]\\[2ex]&=(t_i-t_{i-1})[t_{i-1}, t_i, s, t_j],
    \end{array}
    \]
    we obtain the expression
    \[
    \begin{array}{l}
    D_r(s)=(t_2-t_1)(t_3-t_2)\cdots(t_r-t_{r-1})\times\\[1ex]
    \hskip 5em\left|\begin{array}{cccc}
[t_1, s, t_1]               & [t_1, s, t_2]               & \cdots & [t_1, s, t_r]\\[1ex]
[t_1, t_2, s, t_1]       & [t_1, t_2, s, t_2]       & \cdots & [t_1, t_2, s, t_r]\\[1ex]
        \vdots    & \vdots                    &        & \vdots\\[1ex]
[t_{r-1}, t_r, s, t_1] & [t_{r-1}, t_r, s, t_2 ] & \cdots & [t_{r-1}, t_r, s, t_r]
    \end{array}\right|.
    \end{array}
    \]
    We then subtract the first column from the
    second column, the second column from the third and so forth until the $ (r-1) $th column is subtracted from
    the $ r $th column and obtain
    \[
    \begin{array}{l}
    D_r(s)=(t_2-t_1)^2(t_3-t_2)^2\cdots(t_r-t_{r-1})^2\times\\[1ex]
    \hskip 6em\left|\begin{array}{cccc}
[t_1, s, t_1]                 & [t_1, s, t_1, t_2]                  & \cdots & [t_1, s, t_{r-1}, t_r]\\[1ex]
[t_1, t_2, s, t_1]   & [t_1, t_2, s, t_1, t_2]    & \cdots & [t_1, t_2, s, t_{r-1}, t_r]\\[1ex]
        \vdots      & \vdots                       &        & \vdots\\[1ex]
[t_{r-1}, t_r, s, t_1]& [t_{r-1}, t_r, s, t_1, t_2] & \cdots & [t_{r-1}, t_r, s, t_{r-1}, t_r ]
    \end{array}\right|
    \end{array}
    \]
    since for $ i=1 $ and $ j=2,\dots,r $ the difference
    \[
    \begin{array}{rl}
    [t_1, s, t_j]-[t_1, s, t_{j-1}]&=[t_j, t_1, s]-[t_1, s, t_{j-1}]=(t_j-t_{j-1})[t_j, t_1, s, t_{j-1}]\\[1ex]
    &=(t_j-t_{j-1})[t_1, s, t_{j-1}, t_j],
    \end{array}
    \]
    while for $ i=2,\dots,r $ and $ j=2,\dots,r $ the difference
     \[
    \begin{array}{l}
    [t_{i-1}, t_i, s, t_j]-[t_{i-1}, t_i, s, t_{j-1}]=[t_j, t_{i-1}, t_i, s]-[t_{i-1}, t_i, s, t_{j-1}]\\[1ex]
    =(t_j-t_{j-1})[t_j, t_{i-1}, t_i, s, t_{j-1}]=(t_j-t_{j-1})[t_{i-1}, t_i, s, t_{j-1}, t_j].
    \end{array}
    \]
    Note that the above expression proves (\ref{intermediate determinant identity}) for $ p=1 $
    and any $ r=2,\dots,n. $ In particular
    (\ref{intermediate determinant identity}) is valid for $ r=2. $

    Assume now that $ r\ge 3 $ and (\ref{intermediate determinant identity}) is valid for some $ p\le r-2. $
    We subtract, in the matrix $ M_r(p), $ the $ (p+1) $th row from the $ (p+2) $th row,
    the $ (p+2) $th row from the $ (p+3) $th row until the $ (r-1) $th row is subtracted from the $ r $th row
    and obtain
    \[
    \begin{array}{l}
    \det M_r(p)=(t_{p+2}-t_1)(t_{p+3}-t_2)\cdots(t_r-t_{r-(p+1)})\times\\[2ex]
           \left|\begin{array}{cc}
           M_p & \bigl([t_1,\dots,t_i, s, t_{j-p},\dots,t_j]\bigr)_{i,j}\\[1ex]
           \bigl([t_1,\dots,t_{p+1}, s, t_1,\dots,t_j]\bigr)_{j}
           & \bigl([t_1,\dots,t_{p+1}, s, t_{j-p},\dots,t_j]\bigr)_{j}\\[1ex]
           \bigl([t_{i-(p+1)},\dots,t_i, s, t_1,\dots,t_j]\bigr)_{i,j}
           & \bigl([t_{i-(p+1)},\dots,t_i, s, t_{j-p},\dots,t_j]\bigr)_{i,j}
           \end{array}\right|.
    \end{array}
    \]
    This is so since for $ i=p+2,\dots,r $ and $ j=1,\dots,p $ the difference
    \[
    \begin{array}{l}
    [t_{i-p},\dots,t_i,s,t_1,\dots,t_j]-[t_{i-(p+1)},\dots,t_{i-1},s,t_1,\dots,t_j]\\[1ex]
    =(t_i-t_{i-(p+1)})[t_{i-(p+1)},t_{i-p},\dots,t_i,s,t_1,\dots,t_j]
    \end{array}
    \]
    and similarly
    \[
    \begin{array}{l}
    [t_{i-p},\dots,t_i,s,t_{j-p},\dots,t_j]-[t_{i-(p+1)},\dots,t_{i-1},s,t_{j-p},\dots,t_j]\\[1ex]
    =(t_i-t_{i-(p+1)})[t_{i-(p+1)},t_{i-p},\dots,t_i,s,t_{j-p},\dots,t_j].
    \end{array}
    \]
    Note that row $ p+1 $ was left unchanged.

    We finally subtract, in the above determinant expression, the $ (p+1) $th column from the $ (p+2) $th column,
    the $ (p+2) $th column from the $ (p+3) $th column until the $ (r-1) $th column is subtracted from the
    $ r $th column and obtain
     \[
    \begin{array}{l}
    \det M_r(p)=(t_{p+2}-t_1)^2(t_{p+3}-t_2)^2\cdots(t_r-t_{r-(p+1)})^2\times\\[2ex]
           \left|\begin{array}{cc}
           M_{p+1} & \bigl([t_1,\dots,t_i, s, t_{j-(p+1)},\dots,t_j]\bigr)_{i,j}\\[1ex]
           \bigl([t_{i-(p+1)},\dots,t_i, s, t_1,\dots,t_j]\bigr)_{i,j}
           & \bigl([t_{i-(p+1)},\dots,t_i, s, t_{j-(p+1)},\dots,t_j]\bigr)_{i,j}
           \end{array}\right|
    \end{array}
    \]
    by calculations as above. Our calculations show that
    \[
    \det M_r(p)=(t_{p+2}-t_1)^2(t_{p+3}-t_2)^2\cdots(t_r-t_{r-(p+1)})^2\det M_r(p+1)
    \]
    and consequently
    \[
    D_r(s)=\det M_r(p+1)\cdot\prod_{k=1}^{p+1}\,\,\prod_{l=1}^{r-k}(t_{k+l}-t_l)^2
    \]
    which shows (\ref{intermediate determinant identity}) by induction.

    The second statement of the theorem now follows by choosing $ s=t $ and letting all the numbers
    $ t_1,\dots,t_n $ tend to $ t. $

    If $ f $ is $ n $-convex, then it follows by Kraus' theorem that the matrix $ H(s) $
    is positive semi-definite, thus all the leading determinants (\ref{leading determinants})
    are non-negative. Since the numbers $ t_1,\dots,t_n $ are distinct, it follows
    that also the determinants of the matrices $ M_r $ are non-negative for $ r=1,\dots,n. $
    By choosing $ s=t $ and letting all the numbers $ t_1,\dots,t_n $ tend to $ t, $ we
    derive that the leading principal determinants of the matrix $ K_n(f;t) $ are all non-negative.
    But since each principal submatrix of $ K_n(f;t) $ in this way may be obtained as a leading principal submatrix
    by first making a suitable joint permutation of the rows and columns in $ H(s), $ it follows that the determinants
    of all principal submatrices of $ K_n(f;t) $ are non-negative. Therefore $ K_n(f;t) $
    is indeed positive semi-definite.

    Finally, if the matrix $ K_n(f;t_0) $ is positive definite in some point $ t_0\in I, $
    we use that the entries $ [s,t_i,t_j]_f $ of the matrix $ H(s) $ are continuous functions of $ t_i, t_j $
    and $ s $ to obtain that the matrix $ H(s) $ is positive definite for all $ t_1,\dots,t_n $
    and $ s=t_1,\dots,t_n $ in an open interval $ J $ with $ t_0\in J\subseteq I. $
    The assertion now follows from the characterization by Kraus \cite{kn:kraus:1936}
    of matrix convexity in terms of divided differences.
    \end{proof}

    \subsection{The existence of gaps}

    \begin{proposition}\label{proposition: existence of n-monotone and n-concave polynomials}
    Let $ I $ be a finite interval, and let $ m $ and $ n $ be natural numbers with $ m\ge 2n. $
    There exists an $ n $-concave and $ n $-monotone polynomial $ f_m\colon I\to\mathbf R $ of degree $ m. $
    Likewise there exists an $ n $-convex and $ n $-monotone polynomial $ g_m\colon I\to\mathbf R $ of degree $ m. $
    \end{proposition}

    \begin{proof} We may without loss of generality assume that $ I $ is an open interval and then obtain the
    statement of the proposition for other finite interval types by considering restrictions of a polynomial defined
    on an open interval.
    The interval $ I $ may thus be written on the form $ I=(t_0-c,t_0+c) $ for some $ t_0\in\mathbf R $
    and a positive real number $ c. $
    We introduce the polynomial $ p_m $ of degree $ m $ given by
    \[
    p_m(t)=b_1 t + b_2 t^2 +\cdots+ b_m t^m
    \]
    where
    \[
    b_k=\int_{-1}^0 t^{k-1}\,dt=\frac{(-1)^{k-1}}{k}
    \]
    for $ k=1,\dots,m. $ Thus the $ p $th derivative $ p_m^{(p)}(0)=p!\cdot b_p $ for $ p=1,\dots,2n $
    and consequently
    \[
    M_n(p_m;0)=\left(\frac{p_m^{(i+j-1)}(0)}{(i+j-1)!}\right)_{i,j=1}^n=
    \left(b_{i+j-1}\right)_{i,j=1}^n,
    \]
    where we used the notation in \cite{kn:donoghue:1974,kn:hansen:2004:1}. Similarly
    \[
    K_n(p_m;0)=\left(\frac{p_m^{(i+j)}(0)}{(i+j)!}\right)_{i,j=1}^n=\left(b_{i+j}\right)_{i,j=1}^n.
    \]
    Take a vector $ c=(c_1,\dots,c_n)\in\mathbf C^n, $ then
    \[
    \left(M_n(p_m;0)c\mid c\right)=\sum_{i,j=1}^n b_{i+j-1} c_j\bar c_i=
    \int_{-1}^0\left|\sum_{i=1}^n c_i t^{i-1}\right|^2 dt
    \]
    and
    \[
    \left(K_n(p_m;0)c\mid c\right)=\sum_{i,j=1}^n b_{i+j} c_j\bar c_i=
    \int_{-1}^0 t\left|\sum_{i=1}^n c_i t^{i-1}\right|^2 dt.
    \]
    Since the coefficients in a polynomial are all zero if the polynomial is the zero function, we derive that
    $ M_n(p_m;0) $ is positive definite and $ K_n(p_m;0) $ is negative definite. Hence
    there exists an $ \alpha>0 $ such that $ M_n(p_m;t) $ is positive definite and $ K_n(p_m;t) $ is
    negative definite in the interval $ (-\alpha,\alpha), $ thus $ p_m $ is $ n $-monotone and $ n $-concave
    on $ (-\alpha,\alpha). $ The polynomial
    \[
    f_m(t)=p_m(\alpha c^{-1}(t-t_0))\qquad t\in I
    \]
    then has the desired properties. The second statement is proved by choosing the coefficients
    \[
    b_k=\int_0^1 t^{k-1}\,dt=\frac{1}{k}
    \]
    and then follow the same steps as in the above proof.
    \end{proof}

    \begin{proposition}
    Let $ I $ be an interval, and let $ n\ge 2 $ be a natural number. There are no
    $ n $-convex polynomials of degree $ m $ in $ I $ for $ m=3,\dots,2n-1. $
    \end{proposition}

    \begin{proof}
    If $ f_m $ is an $ n $-convex polynomial of degree $ m $ in $ I $ and $ t_0 $ is an inner point in $ I, $
    then
    \[
    p_m(t)=f_m(t-t_0)
    \]
    is $ n $-convex in a neighborhood of zero and may be written on the form
    \[
    p_m(t)=b_0+b_1 t+\cdots+ b_m t^m
    \]
    where $ b_m\ne 0. $ We calculate the derivatives
    \[
    p_m^{(m-1)}(0)=(m-1)!\, b_{m-1},\quad p_m^{(m)}(0)=m!\, b_m,\quad p_m^{(m+1)}(0)=0.
    \]
    If $ m $ is even and thus of the form $ m=2l $ for some $ l\ge 2, $ then the principal submatrix of $ K_n(p_m; 0) $
    consisting of the rows and columns with numbers $ l-1 $ and $ l+1 $ is given by
    \[
    \left(\begin{array}{cc}
    b_{m-2} &  b_m\\
    b_m     &  0
    \end{array}\right),
    \]
    and it has determinant $ -b_m^2<0. $ If $ m $ is odd and thus of the form $ m=2l+1 $ for some
    $ l\ge 1, $ then the principal
    submatrix of $ K_n(p_m; 0) $ consisting of the rows and columns with numbers $ l $ and $ l+1 $ is given by
    \[
    \left(\begin{array}{cc}
    b_{m-1} &  b_m\\
    b_m     &  0
    \end{array}\right),
    \]
    and this matrix also has determinant $ -b_m^2<0. $ Since $ K_n(p_m, 0) $ is positive semi-definite according
    to Theorem~\ref{determinant condition} we have in both cases a contradiction.
    \end{proof}

    Note that the quadratic polynomial $ t^2 $ is $ n $-convex in any interval for all natural numbers $ n. $

    \begin{corollary}\label{corollary: gaps on finite intervals}
    Let $ I $ be a finite interval, and let $ n $ be a natural number. There exists an $ n $-convex
    function in $ I $ which is not $ (n+1) $-convex in any subinterval of $ I. $
    \end{corollary}

    Note that the function in the corollary may be chosen as
    either an $ n $-monotone increasing or an $ n $-monotone decreasing polynomial of degree $ 2n, $
    and that the possible degrees of any
    polynomials in the gap are limited to $ 2n $ and $ 2n+1. $

    \begin{corollary}
    Let $ I $ be an infinite interval different from the real line. For any natural number $ n $ there is
    an $ n $-convex function in $ I $ which is not $ (n+1) $-convex.
    \end{corollary}

    \begin{proof}
    We may without loss of generality assume that $ I=[0,\infty). $ We may by
    Proposition~\ref{proposition: existence of n-monotone and n-concave polynomials} choose a polynomial
    $ f_n $ of degree $ 2n $ which in the interval $ [0,1) $ is $ n $-monotone and $ n $-concave.
    Possibly by adding a suitable constant we may assume that $ f_n $ is non-negative.

    The transformation
    $ t\to h(t)=t(1+t)^{-1} $ from $ [0,\infty) $ to $ [0,1) $ is operator concave, therefore the function
    \[
    g_n(t)=f_n\bigl(h(t)\bigr)=f_n\left(\frac{t}{1+t}\right)\qquad t\ge 0
    \]
    is $ n $-concave on $ [0,\infty). $ But since the inverse transformation $ t\to h^{-1}(t)=t(1-t)^{-1} $
    from $ [0,1) $ to $ [0,\infty) $ is operator monotone and $ f_n=g_n\circ h^{-1} $ is not
    $ (n+1) $-monotone, we derive that $ g_n $ is not
    $ (n+1) $-monotone.  But a non-negative $ (n+1) $-concave function defined in the interval
    $ [0,\infty) $ is necessarily $ (n+1) $-monotone \cite[Proposition 1.3]{kn:hansen:2003:1}.
    We therefore conclude that $ g_n $ is not $ (n+1) $-concave.
    \end{proof}

    Note that the above proof does not exclude the possibility that $ g_n $ is $ (n+1) $-concave in some
    subinterval of the half-line $ [0,\infty). $

  \section{Local property}

  We say that $ n $-convexity is a local property if an arbitrary function $ f, $ defined in two overlapping open
    intervals $ I_1 $ and $ I_2 $  such that the restrictions of $ f $ to $ I_1 $ and $ I_2 $ are $ n $-convex,
    necessarily is $ n $-convex also in the union $ I_1\cup I_2. $

    We conjecture that $ n $-convexity, like $ n $-monotonocity, is a local property, and we
    prove it for $ n=2. $
   The following representation of divided differences is due to Hermite \cite{kn:hermite:1878}.

  \begin{proposition}\label{Hermite expression}
  Divided differences can be written in the following form
  \[
  \begin{array}{rl}
  [x_0, x_1]_f&=\displaystyle\int_0^1 f'\Bigl((1-t_1)x_0+t_1x_1\Bigr)\,dt_1\\[2ex]
  [x_0, x_1, x_2]_f&=\displaystyle\int_0^1\int_0^{t_1}
  f''\Bigl((1-t_1)x_0+(t_1-t_2)x_1+t_2 x_2\Bigr)\,dt_2\,dt_1\\[2ex]
  &\vdots\\[0pt]
  [x_0, x_1,\cdots, x_n]_f&=\displaystyle\int_0^1\int_0^{t_1}\cdots\int_0^{t_{n-1}}%
  f^{(n)}\Bigl((1-t_1)x_0+(t_1-t_2)x_1+\cdots\\[2ex]
  &\hskip 7em+\,(t_{n-1}-t_n)x_{n-1}+t_nx_n\Bigr)\,dt_n\cdots dt_2\,dt_1
  \end{array}
  \]
  where $ f $ is an $ n $-times continuously differentiable function defined in an open interval $ I, $ and
  $ x_0,x_1,\dots,x_n $ are (not necessarily distinct) points in $ I. $
  \end{proposition}

    \subsection{An inequality for divided differences}

        \begin{proposition}\label{inequality for divided differences}
    Let $ I $ be an open interval and $ n $ a natural number. For a function $ f\in C^{n}(I) $ we
    assume that the $ n $th derivative $ f^{(n)} $ is strictly positive. If in addition the function
    \[
    c(x)=\frac{1}{f^{(n)}(x)^{1/(n+1)}}\qquad t\in I
    \]
    is convex, then the divided difference
    \[
    [x_0,x_1,\dots,x_n]_f\ge\prod_{i=0}^n\: [x_i,x_i,\dots,x_i]_f^{1/(n+1)}
    \]
    for arbitrary $ x_0,x_1,\dots,x_n\in I, $ where the divided differences $ [x_i,x_i,\dots,x_i]_f $ are of
    order $ n. $ If on the other hand the (positive) function $ c(x) $ is
    concave, then the inequality is reversed.
    \end{proposition}

    \begin{proof}
    By using the expression for divided differences given in
    Proposition~\ref{Hermite expression} and the convexity of the function $ c $ we obtain
    \[
    \begin{array}{l}
    [x_0, x_1,\cdots, x_n]_f=
    \displaystyle\int_0^1\int_0^{t_1}\cdots\int_0^{t_{n-1}} f^{(n)}\Bigl((1-t_1)x_0+(t_1-t_2)x_1+\\
    \hskip 11em\cdots+\,(t_{n-1}-t_n)x_{n-1}+t_nx_n\Bigr)
    \,dt_n\cdots dt_2\,dt_1\\[2ex]
    =\displaystyle\int_0^1\int_0^{t_1}\cdots\int_0^{t_{n-1}}
    c\Bigl((1-t_1)x_0+(t_1-t_2)x_1+\\
    \hskip 10em\cdots+\,(t_{n-1}-t_n)x_{n-1}+t_nx_n\Bigr)^{-(n+1)}dt_n\cdots dt_2\,dt_1\\[2ex]
    \ge\displaystyle\int_0^1\int_0^{t_1}\cdots\int_0^{t_{n-1}}\Bigl((1-t_1)c(x_0)+(t_1-t_2)c(x_1)\,+\\
    \hskip 8em\cdots+\,(t_{n-1}-t_n)c(x_{n-1})+t_n c(x_n)\Bigr)^{-(n+1)}dt_n\cdots dt_2\,dt_1.
    \end{array}
    \]
    Next considering the function
    \[
    g(t)=\frac{1}{t}\qquad t>0
    \]
    with $ n $th derivative
    \[
    g^{(n)}(t)=(-1)^n \frac{n!}{t^{n+1}}
    \]
    we may insert this in the above expression to obtain
    \[
    \begin{array}{l}
    \displaystyle [x_0, x_1,\dots, x_n]_f\ge\frac{(-1)^n}{n!}\int_0^1\int_0^{t_1}\cdots\int_0^{t_{n-1}}
    g^{(n)}\Bigl((1-t_1)c(x_0)\,+\\\
    \hskip 3em(t_1-t_2)c(x_1)+\cdots+\,(t_{n-1}-t_n)c(x_{n-1})+t_n c(x_n)\Bigr)
    \,dt_n\cdots dt_2\,dt_1\\[2ex]
    =\displaystyle\frac{(-1)^n}{n!} [c(x_0),c(x_1),\dots,c(x_n)]_g
    \end{array}
    \]
    where we used Proposition~\ref{Hermite expression} once more. Finally, since
    \[
    [t_0,t_1,\dots,t_n]_g=(-1)^n g(t_0)g(t_1)\cdots g(t_n)\qquad t_0,t_1,\dots,t_n>0
    \]
    we obtain
    \[
    \begin{array}{rl}
    [x_0, x_1,\dots, x_n]_f&\ge\displaystyle\frac{1}{n!}\,g(c(x_0))g(c(x_1))\cdots g(c(x_n))\\[2ex]
    &=\displaystyle\frac{1}{n!\, c(x_0)c(x_1)\cdots c(x_n)}\\[2ex]
    &=\displaystyle\frac{1}{n!} f^{(n)}(x_0)^{1/(n+1)} f^{(n)}(x_1)^{1/(n+1)}\cdots f^{(n)}(x_n)^{1/(n+1)}\\[2ex]
    &=\displaystyle\prod_{i=0}^n \left(\frac{f^{(n)}(x_i)}{n!}\right)^{1/(n+1)}
    \end{array}
    \]
    and since $ f^{(n)}(x)=n!\, [x,x,\dots,x]_f $ for $ x\in I $ the statement follows.
    If the function $ c(x) $ is concave, the statement follows by making the appropriate alterations in the
    above proof.
    \end{proof}

    We may use the above Proposition for the exponential function since the function
    \[
    c(x)=\frac{1}{\exp^{(n)}(x)^{1/(n+1)}}=\exp(-x/(n+1))
    \]
    indeed is convex. We therefore obtain
    \[
    \begin{array}{rl}
    [x_0,x_1,\dots,x_n]_{\exp}\ge&\displaystyle\prod_{i=0}^n\: [x_i,x_i,\dots,x_i]_{\exp}^{1/(n+1)}\\[3ex]
    =&\displaystyle\frac{1}{n!}\,\exp\left(\frac{x_0+x_1+\cdots+x_n}{n+1}\right)
    \end{array}
    \]
    for arbitrary real numbers $ x_0,x_1,\dots,x_n. $ This consequence of
    Theorem~\ref{inequality for divided differences} was proved
    in \cite[Appendix 1]{kn:lecouteur:1980} by another method.

  \subsection{Local property for $ 2 $-convex functions}

    \begin{theorem}\label{theorem: local property for 2-convex functions}
    Let $ I $ be an open interval, and take a function $ f\in C^4(I) $ such that $ f''(t)> 0 $ for
    every $ t\in I. $ Then the following assertions are equivalent.

    \begin{list}{(\arabic{std})}{\usecounter{std}}

    \item  $ f $ is 2-convex.

    \item The matrix
    \[
     \left(\begin{array}{cc}
               \displaystyle\frac{f''(t)}{2}  & \displaystyle\frac{f^{(3)}(t)}{6}\\[1.5ex]
               \displaystyle\frac{f^{(3)}(t)}{6} & \displaystyle\frac{f^{(4)}(t)}{24}
               \end{array}\right)
     \]
     is positive semi-definite for every $ t\in I. $

     \item There is a positive concave function $ c $ on $ I $ such that $ f''(t)=c(t)^{-3} $ for every
     $ t\in I. $

     \item The inequality
     \[
     [t_0,t_0,t_0]_f [t_1,t_1,t_1]_f - [t_0,t_1,t_1]_f [t_0,t_0,t_1]_f\ge 0
     \]
     is valid for all $ t_0,t_1\in I. $

     \item The Kraus determinant
     \[
       \left|\begin{array}{cc}
               [t_0,t_0,t_0]_f  & [t_0,t_0,t_1]_f\\[1.5ex]
               [t_0,t_0,t_1]_f  & [t_0,t_1,t_1]_f
               \end{array}\right|\ge 0
     \]
     for all $ t_0,t_1\in I. $

    \end{list}
    \end{theorem}

    \begin{proof} $ (1)\Rightarrow(2) $ is proved by Theorem~\ref{determinant condition}.\\[1ex]
    $ (2)\Rightarrow(3): $ Put $ c(t)=f''(t)^{-1/3} $ for $ t\in I. $ Then $ c $ is a positive
    function and $ f''(t)=c(t)^{-3}. $ By differentiation we obtain $ f^{(3)}(t)=-3 c(t)^{-4}c'(t) $
    and
    \[
    f^{(4)}(t)=12 c(t)^{-5}c'(t)^2-3 c(t)^{-4}c''(t).
    \]
    The determinant
    \[
    \frac{f''(t)}{2}\frac{f^{(4)}(t)}{24}-\frac{f^{(3)}(t)}{36}^2
    \]
    is non-negative by (2), thus inserting the derivatives we obtain
    \[
    \begin{array}{l}
    \displaystyle\frac{c(t)^{-3}}{2}\cdot\frac{12 c(t)^{-5}c'(t)^2-3 c(t)^{-4}c''(t)}{24}
    -\frac{(-3 c(t)^{-4}c'(t))^2}{36}\\[2ex]
    =-\displaystyle\frac{1}{16}c(t)^{-7}c''(t)\ge 0,
    \end{array}
    \]
    hence $ c''(t)\le 0 $ for every $ t\in I $ and $ c $ is concave.\\[1ex]
    $ (3)\Rightarrow(4): $ For $ n=2 $ condition (3) becomes the assumption
    in Proposition \ref{inequality for divided differences}, hence
    \[
    [t_0,t_1,t_2]_f\le [t_0,t_0,t_0]_f^{1/3} [t_1,t_1,t_1]_f^{1/3} [t_2,t_2,t_2]_f^{1/3}
    \]
    for arbitrary $ t_0,t_1,t_2\in I. $ Setting $ t_2=t_0 $ we obtain
    \[
    [t_0,t_0,t_1]_f\le [t_0,t_0,t_0]_f^{2/3} [t_1,t_1,t_1]_f^{1/3}
    \]
    and setting $ t_2=t_1 $ we obtain
    \[
    [t_0,t_1,t_1]_f\le [t_0,t_0,t_0]_f^{1/3} [t_1,t_1,t_1]_f^{2/3},
    \]
    hence the product
     \[
     [t_0,t_1,t_1]_f [t_0,t_0,t_1]_f\le [t_0,t_0,t_0]_f [t_1,t_1,t_1]_f
     \]
     which is condition (4).\\[1ex]
    $ (4)\Rightarrow(5): $
    We introduce a function $ F:I\to\mathbf R $ defined by setting $ F(t_0)=0 $ and
    \[
    \begin{array}{l}
    F(t)=\displaystyle[t_0,t_0,t_0]_f\bigl((t-t_0)f'(t)-f(t)+f(t_0)\bigr)\\[1ex]
    \hskip 12em\displaystyle-\frac{1}{(t-t_0)^2}((t_0-t)f'(t_0)-f(t_0)+f(t))^2.
    \end{array}
    \]
    for $ t\ne t_0. $ Since
    \[
    \frac{1}{(t-t_0)^2}((t_0-t)f'(t_0)-f(t_0)+f(t))^2=(t-t_0)^2[t_0,t_0,t]^2
    \]
    for $ t\ne t_0 $ this defines $ F $ as a differentiable function, and since
    \[
    [t_0,t,t]_f=\frac{1}{(t-t_0)^2}\Bigl((t-t_0)f'(t)-f(t)+f(t_0)\Bigr)
    \]
    we obtain
    \begin{gather}\label{determinant identity}
    [t_0,t_0,t_0]_f[t_0,t,t]_f-[t_0,t_0,t]_f^2=\displaystyle\frac{1}{(t-t_0)^2}F(t)
    \end{gather}
    for $ t\ne t_0. $ We next consider the derivative
    \[
    \begin{array}{rl}
    F'(t)&=\displaystyle [t_0,t_0,t_0]_f\bigl(f'(t)+(t-t_0)f''(t)-f'(t)\bigr)\\[1ex]
    &\hskip 3em\displaystyle +2(t-t_0)^{-3}\bigl((t_0-t)f'(t_0)-f(t_0)+f(t)\bigr)^2\\[1ex]
    &\hskip 3em\displaystyle -2(t-t_0)^{-2}\bigl((t_0-t)f'(t_0)-f(t_0)+f(t))(-f'(t_0)+f'(t)\bigr)\\[2ex]
    &=2(t-t_0)\Bigl([t_0,t_0,t_0]_f[t,t,t]_f-[t_0,t,t]_f[t_0,t_0,t]_f\Bigr).
    \end{array}
    \]
    The assumption $ (4) $ entails that $ F $ has minimum in $ t_0 $ and therefore is non-negative. But this
    is equivalent to $ (5) $ by the identity (\ref{determinant identity}).\\[1ex]
    $ (5)\Rightarrow(1): $ This is the characterization by Kraus \cite{kn:kraus:1936}.
    \end{proof}

      \begin{corollary}
    2-convexity is a local property.
    \end{corollary}

    \begin{proof}
     By applying the regularization procedure described in the introduction we may assume that $ f $
     is infinite many times differentiable. If $ f $ is an affine function there is nothing to prove. If $ f $ is not
     affine we may by Lemma~\ref{lemma: 2-monotone and 2-convex functions} (2) assume that $ f'' $ is strictly
     positive. The statement is now a direct consequence of Theorem~\ref{theorem: local property for 2-convex functions}.
    \end{proof}

    \begin{corollary} Let $ f $ be a twice continuously differentiable function defined in an open interval.
    If the determinant
    \[
     [t_0,t_0,t_0]_f [t_0,t_1, t_1]_f - [t_0,t_0,t_1]_f^2= 0
    \]
    for some $ t_1\ne t_0 $ then also
    \[
    [t_0,t_0,t_0]_f [t_0,t,t]_f - [t_0,t_0,t]_f^2= 0
    \]
    for any $ t $ between $ t_0 $ and $ t_1. $
    \end{corollary}

    \begin{proof}
    We first note that the implications $ (3)\Rightarrow (4) $ and $ (4)\Rightarrow (5) $ in the proof of
    Theorem~\ref{theorem: local property for 2-convex functions} only require the function $ f $ to be
    twice continuously differentiable. The condition in the corollary entails that the function $ F $
    defined in (\ref{determinant identity}) for $ t\ne t_0 $ and with $ F(t_0)=0 $ takes minimum both in $ t_0 $
     and $ t_1 $ and consequently vanishes between the two points.
    \end{proof}

    \section{A fractional transformation}

    Let $ I $ be an open interval and take $ t_0\in I. $ To each function $ f\in C^2(I) $ such that
    $ f'(t)> 0 $ for every $ t\in I, $  Nayak \cite{kn:nayak:2004} considered the following transformation
    \[
    g_{t_0}(t)=-\frac{1}{f(t)-f(t_0)}+\frac{1}{f'(t_0)(t-t_0)}
    \]
    which we write on the form
    \begin{equation}\label{Nayak's transform}
    T(t_0,f)(t)=g_{t_0}(t)=\frac{[t_0,t_0, t]_f}{[t_0,t_0]_f[t_0, t]_f}\qquad t\in I.
    \end{equation}
    The inverse transformation is given by
    \begin{equation}
    f(t)=f(t_0)-\frac{1}{T(t_0,f)(t)-\displaystyle\frac{1}{f'(t_0)(t-t_0)}}
    \end{equation}
    for $ t\ne t_0 $ and $ t\in I. $
    Nayak proved \cite{kn:nayak:2004} the following result:

    \begin{theorem} Let $ n $ be a natural number greater than or equal to two.
    The transform $ T(t_0,f)\in P_n(I) $ for all $ t_0\in I, $ if and only if  $ f\in P_{n+1}(I). $
    \end{theorem}
    Let $ I $ be an open interval and take $ t_0\in I. $ To each
    function $ f\in C^3(I) $ such that
    $ f''(t)> 0 $ for every $ t\in I, $ we consider the following transformation
    \begin{gather}\label{lemma: fractional transformation}
    S(t_0,f)(t)=\frac{[t_0,t_0,t_0,t]_f}{[t_0,t_0,t_0]_f[t_0,t_0,t]_f}\qquad t\in I
    \end{gather}
    with inverse
    \begin{gather}
    f(t)=f(t_0)+f'(t_0)(t-t_0)-\frac{t-t_0}{\displaystyle S(t_0,f)(t)-\frac{1}{[t_0,t_0,t_0]_f (t-t_0)}}
    \end{gather}
    for $ t\ne t_0 $ and $ t\in I. $
    The two transformations are connected in the following way. Consider the function
    $ d_{t_0}\colon I\to\mathbf R $ defined by setting $ d_{t_0}(t)=[t_0,t]_f. $
    Since by a simple calculation
    \[
    [t_0, t]_{d_{t_0}}=[t_0, t_0, t]_f\quad\mbox{and}\quad [t_0, t_0, t]_{d_{t_0}}=[t_0, t_0, t_0, t]_f
    \]
    we obtain
    \begin{equation}\label{connection between the two transforms}
    S(t_0,f)=T(t_0,d_{t_0}).
    \end{equation}
    But Nayak's result is not directly applicable since the function $ d_{t_0} $ depends on $ t_0. $

    \begin{lemma}\label{Sylvester's identity}
    Let $ A=(a_{ij})_{i,j=0,1,\dots,k} $ be a $ (k+1)\times (k+1) $ matrix and consider the $ k\times k $
    matrix $ B=(b_{ij})_{i,j=1,\dots,k} $ defined by setting
    $$
    b_{ij}=a_{00}a_{ij}-a_{i0}a_{0j}\qquad i,j=1,\dots,k.
    $$
    Then the determinant $ \det B=a_{00}^{k-1}\det A. $
    \end{lemma}

    \begin{proof} We may express
    $$
    b_{ij}=\det\left(\begin{array}{cc}
                a_{00} &  a_{0j}\\[1ex]
                a_{i0} &  a_{ij}
                \end{array}\right)\qquad i,j=1,\dots,k
    $$
    and observe that the result follows from Sylvester's determinant identity \cite{kn:akritas:1996}.
    \end{proof}

    \begin{lemma}
    The divided differences of the transform $ S(t_0,f) $ in points $ t_i,t_j\in I $ different from $ t_0 $
    may be written on the form
    \begin{equation}\label{divided differences of S}
    [t_i, t_j]_{S(t_0,f)}=
    \frac{[t_0,t_0]_{d_{t_0}} [t_i, t_j]_{d_{t_0}}-[t_i, t_0]_{d_{t_0}} [t_j,t_0]_{d_{t_0}}}
    {[t_0, t_0]_{d_{t_0}}(d_{t_0}(t_i)-d_{t_0}(t_0))(d_{t_0}(t_j)-d_{t_0}(t_0))},
    \end{equation}
    where as above the function $ d_{t_0}\colon I\to\mathbf R $ is defined by setting $ d_{t_0}(t)=[t_0,t]_f. $
    \end{lemma}

    \begin{proof} By (\ref{connection between the two transforms}) and (\ref{Nayak's transform}) we obtain
    \[
    \begin{array}{rl}
    S(t_0,f)(t)&\displaystyle=\frac{[t_0, t_0, t]_{d_{t_0}}}{[t_0, t_0]_{d_{t_0}} [t_0, t]_{d_{t_0}}}=
    \frac{[t_0,t_0]_{d_{t_0}}-[t_0,t]_{d_{t_0}}}{[t_0,t_0]_{d_{t_0}} (d_{t_0}(t_0)-d_{t_0}(t))}\\[3ex]
    &\displaystyle=\frac{-1}{d_{t_0}(t)-d_{t_0}(t_0)}+\frac{1}{[t_0,t_0]_{d_{t_0}} (t-t_0)}.
    \end{array}
    \]
    We calculate in distinct points $ t_i $ and $ t_j $ (different from $ t_0) $
    \[
    \begin{array}{l}
    (t_i-t_j)[t_i, t_j]_{S(t_0,f)}=S(t_0,f)(t_i)-S(t_0,f)(t_j)=\\[1ex]
    \displaystyle\frac{-1}{d_{t_0}(t_i)-d_{t_0}(t_0)}+\frac{1}{[t_0,t_0]_{d_{t_0}}(t_i-t_0)}
    +\frac{1}{d_{t_0}(t_j)-d_{t_0}(t_0)}-\frac{1}{[t_0,t_0]_{d_{t_0}}(t_j-t_0)}\\[3ex]
     =\displaystyle\frac{d_{t_0}(t_i)-d_{t_0}(t_j)}{(d_{t_0}(t_i)-d_{t_0}(t_0))(d_{t_0}(t_j)-d_{t_0}(t_0))}
    -\frac{t_i-t_j}{[t_0,t_0]_{d_{t_0}} (t_i-t_0)(t_j-t_0)}
    \end{array}
    \]
    from which we obtain (\ref{divided differences of S})
    and then realize that the identity holds in arbitrary points $ t_i,t_j\in I $ different from $ t_0, $
    which is the statement of the lemma.
    \end{proof}

    \begin{theorem}
    Let $ f\in C^3(I) $ where $ I $ is an open interval such that
    $ f''(t)>0 $ for all $ t\in I, $ and let $ n $ be a natural number.
    Then the fractional transform $ S(t_0,f) $ defined in (\ref{lemma: fractional transformation}) is in
    $ P_n(I) $ for all $ t_0\in I, $ if and only if  $ f\in K_{n+1}(I). $
    \end{theorem}

   \begin{proof} Take a $ t_0\in I $ and a sequence $ t_1,\dots,t_n\in I $ of points different from $ t_0. $
   For $ k=1,\dots,n $ we calculate the determinant
   \[
   \begin{array}{l}
   \det\left([t_i, t_j]_{S(t_0,f)}\right)_{i=1}^k\\[2ex]
   =\displaystyle\sum_{\sigma\in S_k}\mbox{Sign}(\sigma)\prod_{i=1}^k
   \left(\frac{[t_0,t_0]_{d_{t_0}}[t_i, t_{\sigma(i)}]_{d_{t_0}}-[t_i, t_0]_{d_{t_0}} [t_{\sigma(i)}, t_0]_{d_{t_0}}}
   {[t_0,t_0]_{d_{t_0}} (d_{t_0}(t_i)-d_{t_0}(t_0))(d_{t_0}(t_{\sigma(i)})-d_{t_0}(t_0))}\right)\\[5ex]
   =\displaystyle\frac{1}{[t_0,t_0]_{d_{t_0}}^{k}}\prod_{j=1}^k\frac{1}{(d_{t_0}(t_j)-d_{t_0}(t_0))^2}
   \sum_{\sigma\in S_k}\mbox{Sign}(\sigma)\\[2ex]
   \hfill\displaystyle\times\prod_{i=1}^k\left([t_0,t_0]_{d_{t_0}} [t_i, t_{\sigma(i)}]_{d_{t_0}}-[t_i, t_0]_{d_{t_0}}
   [t_{\sigma(i)}, t_0]_{d_{t_0}}\right)\\[3ex]
   =\displaystyle\frac{\det B}{[t_0, t_0]_{d_{t_0}}^{k}}\prod_{j=1}^k\frac{1}{(d_{t_0}(t_j)-d_{t_0}(t_0))^2}
   \end{array}
   \]
   where $ B=(b_{ij})_{i,j=1}^k $ and
   \[
   b_{ij}=[t_0, t_0]_{d_{t_0}}[t_i, t_j]_{d_{t_0}}-[t_i, t_0]_{d_{t_0}} [t_0, t_j]_{d_{t_0}}
   \]
   for $ i,j=1,\dots,k. $ If we set $ A=(a_{ij})_{i,j=0}^k $ where
   $ a_{ij}=[t_i,t_j]_{d_{t_0}} $ we may write
   \[
   b_{ij}=a_{00}a_{ij}-a_{i 0} a_{0 j}\qquad i,j=1,\dots,k.
   \]
   By applying Lemma \ref{Sylvester's identity} we obtain
   \[
   \det B=a_{00}^{k-1}\det A=[t_0,t_0]_{d_{t_0}}^{k-1}\det A
   \]
   and by using the identity
   \[
   \det A=\det\left([t_i,t_j]_{d_{t_0}}\right)_{i=0}^k=\det\left([t_0,t_i,t_j]_f\right)_{i,j=0}^k
   \]
   we may write   
   \[
   \det\left([t_i, t_j]_{S(t_0,f)}\right)_{i,j=1}^k
   =\frac{\det\left([t_0,t_i,t_j]_f\right)_{i,j=0}^k}{[t_0,t_0]_{d_{t_0}}}
   \prod_{j=1}^k\frac{1}{(d_{t_0}(t_j)-d_{t_0}(t_0))^2}.
   \]
   Suppose $ f $ is $ (n+1) $-convex.  
   We first obtain $ \det\left([t_i, t_j]_{S(t_0,f)}\right)_{i,j=1}^k\ge 0 $
   for $ k=1,\dots,n. $ By considering permutations of $ t_1,\dots,t_n $ we realize 
   that the determinants of all the principal submatrices of each order $ k $ of the Pick matrix 
   $ \left([t_i, t_j]_{S(t_0,f)}\right)_{i,j=1}^n $ are non-negative, hence
   the Pick matrix itself is positive semi-definite. By continuity we obtain that the Pick matrix
   is positive semi-definite for
   arbitrary sequences $ t_1,\dots,t_n $ in $ I, $ hence $ S(t_0,f) $ is $ n $-monotone.

   If on the other hand $ S(t_0,f) $ is $ n $-monotone, we realize that
   \[
   \det\left([t_0,t_i,t_j]_f\right)_{i,j=0}^k\ge 0
   \]
   for $ k=1,\dots,n $ and since for $ k=0 $ the entry $ [t_0,t_0,t_0]_f>0, $ we obtain that the leading
   determinants of the matrix $ \left([t_0,t_i,t_j]_f\right)_{i,j=0}^k $ are non-negative. By considering
   permutations of $ t_1,\dots,t_n $ we realize that the determinants of all the principal submatrices 
   of each order $ k $ of $ \left([t_0,t_i,t_j]_f\right)_{i,j=0}^n $ are non-negative,
   hence the matrix  is positive semi-definite. By continuity we finally realize that the Kraus matrix in equation 
   (\ref{Kraus matrix}) for $ f $ calculated in arbitrary $ n+1 $ points $ t_0,t_1,\dots,t_n\in I $ is
   positive semi-definite, hence $ f $ is $ (n+1) $-convex.
   \end{proof}

{\footnotesize



      \vfill

      \noindent Frank Hansen: Department of Economics, University
       of Copenhagen, Stu\-die\-straede 6, DK-1455 Copenhagen K, Denmark.\\[1ex]

       \noindent Jun Tomiyama: Department of Mathematics and Physics. Japan Women's University.
       Mejirodai Bunkyo-ku, Tokyo, Japan.

       }

      \end{document}